\documentstyle[11pt]{article}
\begin{document}
\makeatletter \@addtoreset{equation}{section} \makeatother
\renewcommand{\theequation}{\thesection.\arabic{equation}}
\baselineskip 15pt

\title{\bf Book Review: \\ ``Geometry of Quantum States'' \\
by Ingemar Bengtsson and Karol \.Zyczkowski \\ {\small\bf Cambridge University Press
(2006)}}
\author{Paul B. Slater\footnote{e-mail: slater@kitp.ucsb.edu}\\
{\small ISBER,} \\
{\small University of California, Santa Barbara, 93106.}}
\date{}
\maketitle

The authors, distinguished mathematical physicists, have written a markedly distinctive, dedicatedly pedagogical, suitably rigorous text, designed, in part, for advanced undergraduates familiar with the principles of quantum mechanics. The book, pleasing in character and enthusiastic in tone, has many stimulating diagrams and tables, as well as problem sets (with hints and answers supplied at the end). 
The diverse topics covered --- conveniently all assembled here --- reflect the geometrically-oriented, fundamental quantum-information-theoretic interests and expertise of the two authors. 
(Several of the areas surveyed are among the many discussed in another recent copiously illustrated book, Roger Penrose's The road to reality MR2116746 (2005k:83002), which has a broader, less specialized audience as its overall target.)

In the first chapter (``Convexity, colours and statistics''), convexity is discussed using, most notably, the three-dimensional example of color-mixing --- making use of chromaticity diagrams and MacAdam ellipses. The second chapter (``Geometry of probability distributions'') deals with the geometry of probability distributions, setting the stage for the later quantum extension of these concepts. (The reviewer has examined a certain area of overlap between the classical and quantum treatments, in this regard,  J. Math. Phys.  47  (2006),  no. 2, 022104, MR2208138 (2007f:94015).) The next two chapters (``Much ado about Spheres'' and ``Complex Projective Spaces'') are quite technical and sophisticated, discussing the Hopf fibration and Segre embedding, among other topics. In the fifth chapter (``Outline of quantum mechanics''), the authors pause to review some important finite-dimensional aspects of quantum mechanics. The level of sophistication continues high in the next nine chapters: (6) ``Coherent states and group actions''; (7) ``The stellar representation'' (a subject of obvious strong interest to Penrose); (8) ``The space of density matrices'', the one chapter having a summary section; (9) ``Purification of mixed quantum states''; (10) ``Quantum operations''; (11) ``Duality: maps versus states''; (12) ``Density matrices and entropies''; (13) ``Distinguishability measures''; and (14)  ``Monotone metrics and measures''. 

Of particular note, for its involvement with the current frontier of research, is the final (fifteenth) chapter, ``Quantum entanglement'', some fifty pages in length. The authors judiciously choose to limit their coverage here to the case of bipartite entanglement, since the area of multipartite entanglement is still experiencing rapid development --- for example, the uncovering of relations, pertaining to Cayley's hyperdeterminant, between entanglement and strings and black holes (R. Kallosh and A. Linde, Phys. Rev. D (3)  73  (2006),  no. 10, 104033 MR2224730 (2007d:83149);  
Peter L{\"e}vay,  Phys. Rev. D (3)  74  (2006),  no. 2, 024030 MR2249975; Phys. Rev. D  75 (2007), 024024]). Apparently obtained too recently for inclusion in Chapter 15 is a certain (bipartite) result of S. Szarek together with Bengtsson and {\.Z}yczkowski. They established that the set of positive-partial-transpose (or, equivalently, the set of separable in the two smallest possible [qubit-qubit and qubit-qutrit] cases) states is ``pyramid decomposable'' and, hence, is a body of constant height (J. Phys. A  39  (2006),  no. 5, L119--L126 MR2200422 (2006i:81029) ).  The authors do generously refer to and elaborate upon various important papers of W. K. Wootters, but perhaps they might have expanded upon their remark (p. 152), 
as to the ``rather peculiar measures'' taken in Wootters' (much-cited) paper, ``A Wigner-function formulation of finite-state quantum mechanics'',  Ann. Physics  176  (1987),  no. 1, 1--21. (MR0893477 (89d:81047)).

The Epilogue consists of four appendices, of particular note being Appendix 3, ``Geometry: do it yourself'', in which some cutting and gluing exercises are suggested. The other appendices are: ``Basic notions of differential geometry'' (1); ``Basic notions of group theory'' (2); and ``Hints and answers to the exercises'' (4). The impressive References list consists of more than 600 items, most of recent origin.

The stress throughout the text is upon finite-dimensional state spaces, discrete variables and static scenarios. Passing mention only (p. 234) is made of identical (bosonic, fermionic, anyonic) particle entangled systems.

The (subject) index to the book is somewhat slim in the number of main headings 
employed, but certainly rich in subheadings. There is no name index nor (not strongly needed) list of notation --- as is perhaps standard for a textbook. There is a typographical glitch at the end of the line preceding eq. (9.42). The English usage is at a high, eminently readable level throughout, and I detected only a few very minor slips in this regard.

Another substantial, praiseworthy 2006 (exercise-including) volume (Masahito Hayashi [Quantum information, Springer-Verlag, 2006; MR2228302 (2007b:81050)]) --- having a far smaller number of visual aids and lesser stress on geometry {\it per se}  --- has a significant amount of subject matter in common with the 
Bengtsson-{\.Z}yczkowski (BZ) monograph. The Hayashi monograph (H)  is less discursive and more formal than BZ, and directed more for graduate students than undergraduates. The formality of H is, however, considerably offset by interspersed discussions and the very detailed, interesting Historical Notes at the ends of Chapters 3-10, as well as the Summaries at the outset of each chapter. Quantum mechanical issues {\it per se} are not discussed until Chap. 5 [p. 131] of BZ, while they are from the beginning of H. The  References list in H of over 400 items does not include any of the numerous interesting items of {\.Z}yczkowski and his co-authors (in particular, those with H.-J. Sommers) --- such as the computation of the volumes and boundary hyperareas, in terms of the Hilbert-Schmidt and Bures metrics, of the spaces of $n \times n$ density matrices (cf. Attila Andai, J. Phys. A 39 (2006), no. 44, 13641--13657 MR2265929 (2007f:81028)) --- treated, at length,  in BZ (secs. 14.3 and 14.4). Reciprocally, the BZ list has no Hayahsi-authored papers.
(This absence of mutual citations surely reflects more upon the vastness of the quantum-information-theoretic literature than any neglect upon the parts of the authors. H draws heavily upon the Japanese literature, and BZ makes profitable use, in particular, of the rich body of Polish-authored contributions.)  
It appears that H makes no reference to the use of the Bayesian paradigm in quantum information, while BZ does limitedly, referring (p. 325) to ``Bayes 'decision rule'' and (p. 55) to {\it Jeffrey's prior} 
(cf. C. Krattenthaler  and P. B. Slater, IEEE Trans. Inform. Theory 46 (2000), no. 3, 801--819.   MR1763463 (2001e:94010) for a certain geometric Bayesian setting, and P. B. Slater, Phys. Lett. A 247 (1998), no. 1-2, 1--8. MR1650497 (99i:81025) for the question of selecting noninformative priors).

Both BZ and H only make passing mention of quantum computation, a central topic of the majority of recent monographs reporting quantum information processing advances (e. g., M. A. Nielsen and I. L. Chuang [ Quantum computation and quantum information, Cambridge Univ. Press, Cambridge, 2000; MR1796805 (2003j:81038)]; P. Lambropoulos and D. Petrosyan [ Fundamentals of quantum optics and quantum information, Springer, Boston, 2007; see also for {\it topological} quantum computation,  M. H. Freedman, A. Kitaev, M. Larsen and Z. Wang, Bull. Amer. Math. Soc. (N.S.)  40  (2003),  no. 1, 31--38 MR1943131 (2003m:57065)]). (R. Jozsa and N. Linden,
R. Soc. Lond. Proc. Ser. A Math. Phys. Eng. Sci.  459  (2003),  no. 2036, 2011--2032 (MR1993666 (2004f:81055)), among others, have analyzed the apparently quite subtle role of entanglement {\it per se} in quantum computation.)

If widely adopted as a text --- as seems certainly appropriate --- BZ (dense and encyclopedic in scope) should serve to stimulate
substantial numbers of students to investigate research topics in the rapidly growing, challenging domain of quantum information theory. 
An instructor might, depending upon the goals and backgrounds of the students, choose to supplement the book with some quantum computing literature.

\end{document}